\documentclass{article}
\usepackage[utf8]{inputenc}
\usepackage{amsmath, amsfonts, amssymb}

\usepackage[utf8]{inputenc}
\usepackage[english, russian]{babel}
\usepackage{todonotes}
\usepackage{comment}
\usepackage{enumitem}
\usepackage[left=2cm, right=3cm, top=2cm, bottom=2cm, bindingoffset=0cm, marginparwidth=3cm]{geometry}

\usepackage{hyperref}
\usepackage{float}

\usepackage{amsthm}

\usepackage{tikz}
\usetikzlibrary{arrows,calc}
\tikzset{
>=stealth',
}

\theoremstyle{definition}

\theoremstyle{remark}

\newtheorem*{rem*}{Remark}

\addto\captionsrussian{}
\addto\captionsenglish{\renewcommand{\refname}{Bibliography}}

\makeatletter
\def\blfootnote{\gdef\@thefnmark{}\@footnotetext}
\makeatother

\title{The problem of stochastic system prediction in gait biomechanics applications}
\author{S. S. Gavryushin, I. A. Meshchihin, S. S. Minkov}
\date{August 2026}

\begin{document}
\maketitle
\begin{abstract}
This paper analyzes the stochastic dynamics of a two-link pendulum that models lower limb segments (the thigh and shank) during the swing phase of the gait cycle. We describe the model using a system of non-linear Itô stochastic differential equations in a four-dimensional phase space under both external and internal random perturbations. To analyze how the phase-state probability density evolves, we formulate a boundary value problem for the four-dimensional Fokker–Planck equation with zero-flux boundary conditions. We then develop a numerical solution algorithm based on operator splitting and the conservative Chang–Cooper scheme. Finally, we perform a spectral analysis of the discrete Kolmogorov operator, compute the stationary invariant measure, and estimate diffusion parameters using Neural Stochastic Differential Equations (Neural SDEs).

    \vspace{0.5em} 
\end{abstract}

\section{Introduction}
When designing control systems, particularly electronically controlled prostheses, it is essential to analyze how the dynamic system responds to external actions \cite{winter2009}. These actions often exhibit a stochastic nature \cite{gardiner1985}. To analyze a system under these conditions, several critical factors must be combined: a limit state defined by the probability of a specific event (such as the probability of tripping), the inherent non-linearity of the dynamic system model, and the stochastic nature of the system parameters, loading conditions, and initial states. Accounting for all these factors simultaneously requires formulating and solving the dynamic problem as a system of stochastic differential equations \cite{risken1989}.

This paper aims to model the dynamics of a prosthesis during the swing phase using a system of non-linear It\^{o} stochastic differential equations (SDEs) \cite{gardiner1985, oksendal2003}. We identify the drift and diffusion coefficients of this model directly from high-resolution empirical telemetry data. To ensure accuracy, we apply regularization procedures based on prior physical constraints derived from the two-link pendulum equations.

\section{Mathematical Model of the Two-Link Pendulum}

We consider the "thigh-shank" biomechanical system during the swing phase as a two-link pendulum \cite{winter2009, zagrevsky2007} moving in the sagittal plane under the action of control and dissipative torques.

\begin{figure}[htbp]
    \centering
    \includegraphics[width=0.6\textwidth]{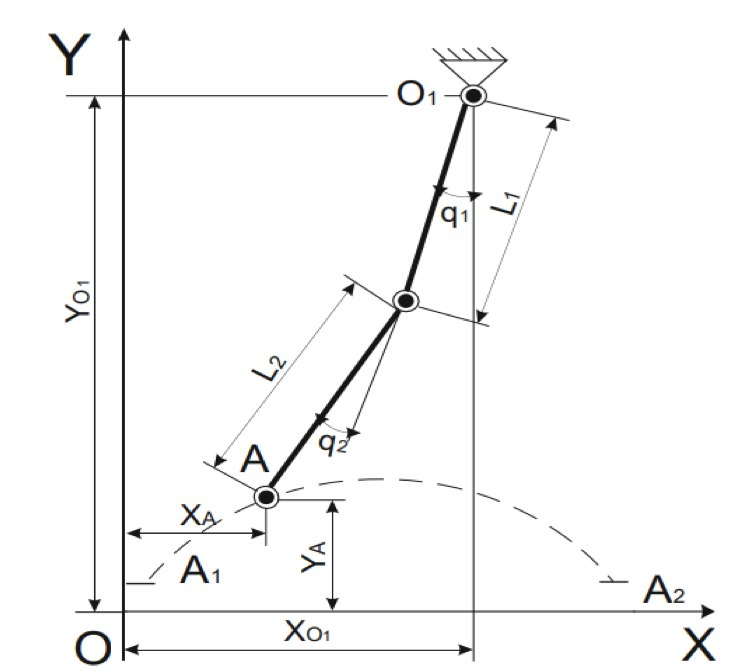}
    \caption{Kinematics of the swing phase.}
    \label{fig:kin}
\end{figure}

Let us define the parameters. The two links consist of the thigh (length $L_1$, mass $m_1$, moment of inertia $I_1$) and the shank (length $L_2$, mass $m_2$, moment of inertia $I_2$). The generalized coordinates are chosen as $\theta_1$, the hip joint angle, and $\theta_2$, the knee joint angle. The control parameter is denoted by $u$, which represents the resistance torque in the knee module (governed by the spool valve position). 

Thus, the deterministic system of equations of motion is given by:
\begin{equation}
    M(\theta)\ddot{\theta} + C(\theta, \dot{\theta})\dot{\theta} + G(\theta) = Q(u)\dot \theta,
    \label{eq:deterministic_motion}
\end{equation}
where $M(\theta)$ is the mass matrix, $C(\theta, \dot{\theta})$ is the matrix of Coriolis and centrifugal forces, $G(\theta)$ is the gravity vector, and $Q(u) \dot \theta$ is the generalized force vector accounting for control and dissipative torques.

Accounting for random perturbations, the system is subjected to a stochastic vector process $\xi(t)$. Transforming this into the It\^{o} form yields:
\begin{equation}
    dx(t) = a(x, u)dt + b(x, u)dW_t,
    \label{eq:ito_sde}
\end{equation}
where $x = (\theta_1, \theta_2, \omega_1, \omega_2)^T$ is the four-dimensional state vector (with angular velocities $\omega_i = \dot{\theta}_i$), $a(x, u)$ is the drift coefficient vector, $b(x, u)$ is the diffusion coefficient matrix, and $W_t$ denotes the standard Wiener process.

\section{Identification of Drift and Diffusion from Data}

In practical applications, the precise analytical identification of the components of matrices $M,C,G$ and $a,b$ is hindered by significant uncertainty in anthropometric characteristics. Consequently, various data-driven approaches are employed to approximate the right-hand sides of the stochastic differential equations (SDEs).

Given the low dimensionality of the state space, it is computationally rational to estimate acceleration values on a grid of nodes as a function of coordinates using kernel regression methods. This explicit approximation of acceleration as a function of phase coordinates allows the system to account for global behavior 'far away' from the observed data density while satisfying kinematic constraints. Figure 2 illustrates the acceleration dependence for a single-link approximation of the hip joint angle (under the assumption that the second link is negligible). Far from the training data domain, the system is approximated purely by a simple pendulum model.

\begin{figure}[htbp]
    \centering
     \includegraphics[width=0.8\textwidth]{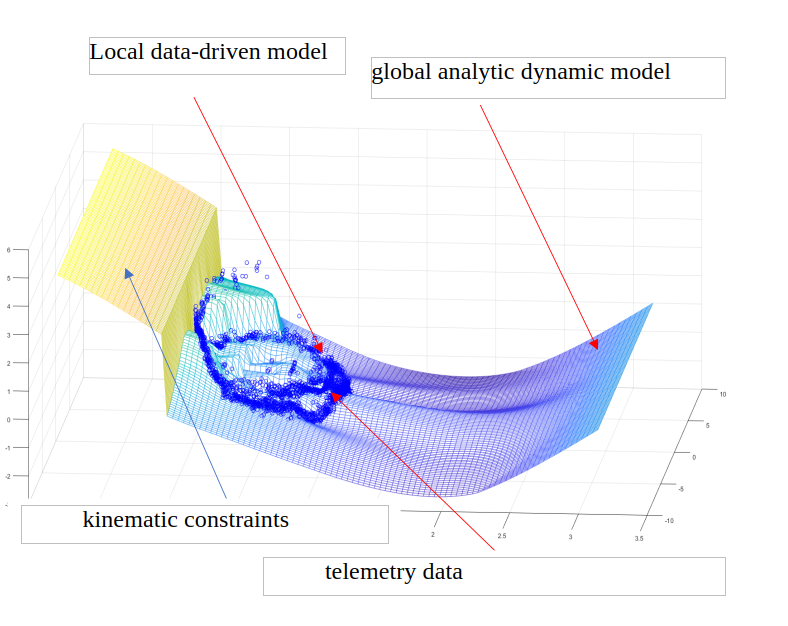}
    \caption{Acceleration dependence for the single-link approximation of the hip joint angle. The model reverts to a pure pendulum behavior outside the experimental data range.}
    \label{fig:acceleration_fit}
\end{figure}

In the multidimensional case, where kernel estimation suffers from the curse of dimensionality, the Neural SDE approach \cite{kidger2020} proves highly effective. Within this framework, the deterministic angular accelerations (drift) and the intensities of random fluctuations (diffusion) are modeled as parametric functions using deep multilayer neural networks. Throughout the rest of this work, we operate strictly within this multidimensional, neural SDE framework unless stated otherwise.

The parameterization of the weight coefficients is performed based on the numerical integration of the stochastic system via the Euler–Maruyama scheme \cite{oksendal2003}. Specifically, the neural network parameters are optimized such that, following the Euler–Maruyama integration step, the resulting path distribution minimizes the expected loss of the objective functional. The optimization objective functional is formulated within the framework of Physics-Informed Neural Networks (PINN) \cite{raissi2019} and is defined as a linear combination of the mean squared error (MSE) of the empirical phase trajectories $L_{data}$ and a penalty function penalizing violations of prior anatomical and kinematic constraints $L_{physics}$:
\begin{equation}
    L_{total} = \lambda_{data} L_{data} + \lambda_{physics} L_{physics}.
    \label{eq:pinn_loss}
\end{equation}

The penalty term $L_{physics}$ imposes a quadratic cost on outlying points if the calculated trajectories exceed the allowable joint angle ranges or maximum physiological angular velocities of the links. This constraint guarantees the internal stability and physical interpretability of the synthesized neural network stochastic model \cite{raissi2019}.

To train the model, a dataset of gait telemetry was experimentally recorded under various speeds and spool valve position setpoints (which dictate how the resistance torque depends on the angular velocity in the knee joint). The training process yields two parametric functions. Given the angles and angular velocities of both the knee and hip joints along with the valve position setpoint, these functions map the inputs to the estimated angular accelerations (drift) and their corresponding variances (diffusion) for both joints.

\begin{figure}[htbp]
    \centering
    \includegraphics[width=0.7\textwidth]{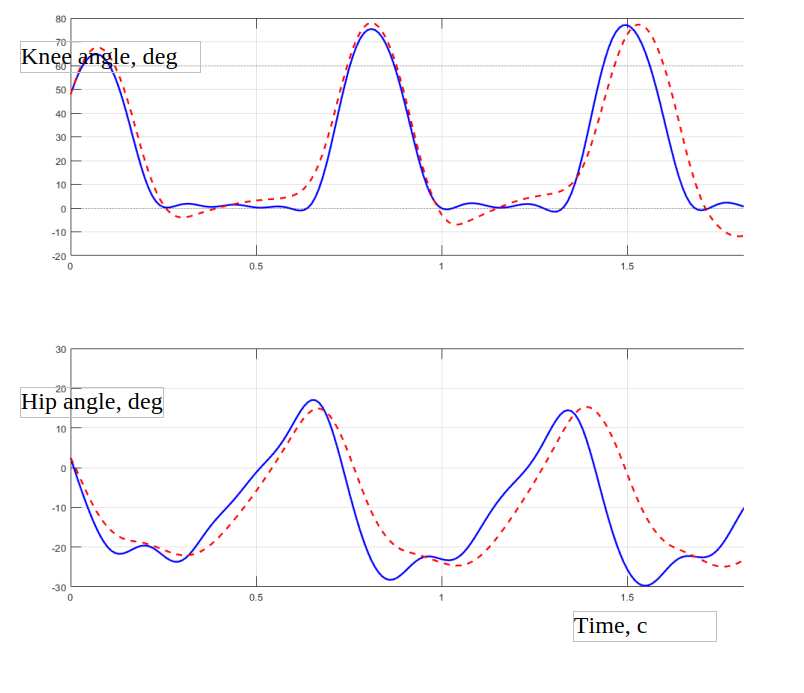}
    \caption{Comparison between experimental telemetry data (solid line) and the model forecast (dashed line).}
    \label{fig:telemetry_vs_prediction}
\end{figure}

As illustrated in Figure 3, the forecast closely reproduces the experimental telemetry data. Although the full telemetry dataset is utilized, the gait cycle phases are explicitly segmented; consequently, during the stance phase, the loss function is evaluated with very small penalty coefficients $\lambda$.

\section{The Fokker–Planck Equation and the Boundary Value Problem}

The drift and diffusion operators obtained via the aforementioned approach are subsequently utilized to solve the Fokker–Planck equation. Unlike a mechanical two-link pendulum, the true dimensionality of the phase space for human gait is unknown a priori. To test the validity of the chosen dimensionality, one can analyze the autocorrelation function of the residual noise (i.e., the remaining signal obtained after subtracting the estimated drift from the empirical telemetry data). If the true system dimensionality exceeds the proposed four dimensions, or if the underlying dynamics cannot be accurately captured by differential equations (requiring instead integro-differential or time-delay equations), the autocorrelation function will deviate from that of pure white noise. In practice, this deviation can be evaluated using various spectral distance metrics.

We consider two approximations. 

Assuming the noise is sufficiently close to a white noise process (which is justified by standard modeling conventions), the classical Fokker–Planck equation can be applied:
\begin{equation}
    \frac{\partial p}{\partial t} = \sum_i \frac{\partial}{\partial x_i} \big( p \, a_i(x) \big) + \sum_i \sum_j \frac{\partial^2}{\partial x_i \partial x_j} \big( p \, D(x_i, x_j) \big).
    \label{eq:fokker_planck}
\end{equation}

The \textbf{drift vector} takes the form $a(x) = (a_1, a_2, a_3, a_4)^T$. Its first two components are strictly determined by kinematic constraints (mapping velocities directly from coordinates). The angular accelerations are obtained by solving the system of Lagrange equations of motion with respect to the highest-order derivatives.

The \textbf{diffusion matrix}: Because random perturbations act on the system exclusively as generalized forces and torques within the system dynamics, the diffusion matrix is degenerate. The coordinate subspaces corresponding to the angles are diffusion-free ($D_{ii} = 0$ for $i=1,2$), whereas the non-zero coefficients are localized entirely within the angular velocity block (e.g., the direct variances $D_{33}$, $D_{44}$ and the cross-diffusion term $D_{34}$).

The second approximation accounts for the case where the noise autocorrelation function is more accurately fitted by the following exponential form \cite{kac1974, stratovich1961}:
\begin{equation}
    K_{\xi}(\tau) = \langle \xi(t) \xi(t + \tau) \rangle = \frac{D}{\tau_c} \exp \left( -\frac{|\tau|}{\tau_c} \right),
    \label{eq:colored_noise_autocorr}
\end{equation}
where $D$ represents the intensity of the random perturbation, and $\tau_c$ denotes the characteristic correlation time (memory) of the process. In the limit as $\tau_c \to 0$, the autocorrelation function converges to that of a delta-correlated pure white noise process.

Under these conditions, the evolution of the phase-state probability density is governed not by the parabolic Fokker–Planck equation, but rather by a generalized hyperbolic telegraph-type equation \cite{kac1974, stratovich1961}:
\begin{equation}
    \tau_c \frac{\partial^2 p}{\partial t^2} + \frac{\partial p}{\partial t} = -\sum_{i=1}^4 \frac{\partial}{\partial x_i} \big( a_i(x) p \big) + \sum_{i,j=3}^4 \frac{\partial^2}{\partial x_i \partial x_j} \big( D_{ij}(x) p \big).
    \label{eq:telegraph_pde}
\end{equation}

Physically, adopting the telegraph equation accounts for the finite propagation speed of probability fronts in the phase space, which scales as $v \sim \sqrt{D / \tau_c}$. Throughout the remainder of this paper, we assume the noise is close to a white noise process (the first case). However, we strictly select numerical algorithms that can be seamlessly generalized to the hyperbolic case (the telegraph equation).

The dimensionality of the problem is fixed at four, with the computational domain defined as:
\begin{equation}
    \Omega = [-\pi, \pi] \times [-\pi, \pi] \times [\omega_{h,\min}, \omega_{h,\max}] \times [\omega_{k,\min}, \omega_{k,\max}].
    \label{eq:domain}
\end{equation}

The boundary conditions are formulated as follows: periodic boundary conditions are applied along the angular coordinates to reflect the cyclic nature of the gait cycle, whereas zero-flux (no-flux) boundary conditions for the probability density are imposed at the velocity boundaries at infinity (or, alternatively, symmetry conditions are utilized). Note that the angular velocity boundaries $\omega_{\min}$ and $\omega_{\max}$ are chosen to be sufficiently large to fully encompass the entire physiological range of human gait; consequently, the probability density asymptotically approaches zero at these boundaries, justifying the zero-flux condition.

\subsection{Numerical Solution}
To solve the four-dimensional Fokker–Planck equation numerically, we implement a hybrid finite-difference solver based on the directional operator splitting method \cite{marchuk1989}. This approach overcomes the "curse of dimensionality" inherent to multidimensional grid-based methods by sequentially applying optimized one-dimensional computational schemes for each phase variable at every time step. 

Furthermore, under this directional splitting framework, the degeneracy of the Fokker–Planck diffusion matrix transforms from a catastrophic bottleneck into a significant computational advantage. The integration algorithm is structured into the following distinct stages:

\begin{itemize}
    \item \textbf{High-resolution transport along angular coordinates:} For the spatial displacement variables, we apply a periodic Total Variation Diminishing (TVD) scheme \cite{harten1983}. This scheme utilizes a Lax–Wendroff anti-diffusive correction combined with the van Leer flux limiter \cite{vanleer1979}. It accurately models convective transport (advection) without generating unphysical oscillations, while strictly preserving the periodic boundary conditions over the full rotation range $[-\pi, \pi]$. To maintain numerical stability at high velocities, an automatic time-step sub-cycling mechanism is implemented, ensuring compliance with the Courant–Friedrichs–Lewy condition ($\text{CFL} < 0.8$).
    
    \item \textbf{Diffusion and drift along velocities:} The evolution of the system along the angular velocity coordinates is computed using the Chang–Cooper scheme \cite{chang1970}. This scheme is selected for its ability to accurately balance convection and diffusion while naturally enforcing the zero-flux boundary conditions at the edges of the computational domain. Each splitting step yields a system of linear algebraic equations with a tri-diagonal matrix, which is efficiently solved using the Thomas algorithm (Tridiagonal Matrix Algorithm, TDMA) \cite{marchuk1989}.
    
    \item \textbf{Cross-diffusion treatment:} The mutual interaction between the velocity components is accounted for in a separate splitting step via an explicit finite-difference approximation of the second-order mixed derivatives.
\end{itemize}

The implemented hybrid solver architecture guarantees the mass-conservative property of the numerical scheme, ensuring strict renormalization and total probability conservation throughout the entire forecasting horizon. To enhance computational efficiency, the core matrix operations are fully vectorized across independent dimensions.

For the angular coordinates, a pure advection equation is solved:
\begin{equation}
    \frac{\partial p}{\partial t} + \omega \frac{\partial p}{\partial \theta} = 0,
    \label{eq:pure_advection}
\end{equation}
where the advection velocity $\omega$ is sampled from the corresponding node of the velocity grid.

To prevent unphysical oscillations and numerical diffusion (front blurring), we implement a Total Variation Diminishing (TVD) scheme. This scheme blends the first-order upwind and second-order Lax–Wendroff schemes using a flux limiter:
\begin{equation}
\end{equation}

\textbf{Sub-cycling (CFL Protection):} The algorithm automatically subdivides the global time step $\Delta t$ into $n$ smaller sub-steps to ensure that the local Courant number satisfies the strict stability condition:
\begin{equation}
    C = \frac{\omega \Delta t_{sub}}{\Delta \theta} < 0.8.
    \label{eq:cfl_condition}
\end{equation}

\textbf{Flux calculation:} The total numerical flux across the cell face is given by:
$$ 
F_{j+1/2}=F_{low}+\phi (r_{j})F_{corr}
$$

\textbf{First-order upwind flux:}
$
    F_{low} = \max(0, u) p_j + \min(0, u) p_{j+1}.
    $

\textbf{Anti-diffusive correction (Lax-Wendroff):}
\begin{equation}
    F_{corr} = \frac{1}{2} |u| \left( 1 - |c_{\text{sub}}| \right) (p_{j+1} - p_j) \operatorname{sign}(u)
\end{equation}

\textbf{Van Leer slope limiter:}
\begin{equation*}
    \phi(r) = \frac{r + |r|}{1 + |r| + \epsilon}.
\end{equation*}

\textbf{Periodic boundary conditions:} Spatial shift indices are mapped onto a ring using the modulo operator:
\begin{equation}
    j_{p1} = \bmod(j, N) + 1, \quad j_{m1} = \bmod(j - 2, N) + 1
\end{equation}
Due to this, the support of the probability measure remains between $-\pi$ and $\pi$ at any given time.

For velocities, physical advection (acceleration) and diffusion (noise) are computed simultaneously:
\begin{equation}
    \frac{\partial p}{\partial t} + \frac{\partial F_{\text{CC}}}{\partial \omega} = 0
\end{equation}
Here, it is critically important to maintain the non-negativity of the density and the exact balance in the steady state (Maxwell–Boltzmann distribution). To achieve this, the Chang–Cooper scheme is implemented in the code.

The local Peclet number (balance parameter \(w\)) at the cell faces is calculated as the ratio of drift \(b\) to diffusion \(D\):
\begin{equation}
    w_{j+1/2} = \frac{b_{j+1/2} \Delta \omega}{D_{j+1/2}}
\end{equation}
where
\begin{itemize}
    \item $w \gg 1$ --- deterministic drift dominates (the system is ``predictable'').
    \item $w \sim 1$ --- balance of drift and diffusion (the zone of maximum uncertainty).
    \item $w \ll 1$  --- random diffusion dominates (the system is ``stochastic'' or ``random'').
\end{itemize}

The Chang–Cooper weight function regularizes the scheme between the central difference and upwind schemes:
$$
f_{\text{CC}}(w) = \frac{w}{1 - e^{-w}}.
$$
The numerical flux is written in the following form:
\begin{equation}
    F_{j+1/2} = \frac{D_{j+1/2}}{\Delta \omega} \left[ f_{\text{CC}}(w_{j+1/2}) p_j - f_{\text{CC}}(-w_{j+1/2}) p_{j+1} \right]
\end{equation}

The time step is implemented using a fully implicit Euler scheme, which leads to a tridiagonal system of equations for each one-dimensional slice:

\begin{equation}
    A_j p_{j-1}^{n+1} + B_j p_j^{n+1} + C_j p_{j+1}^{n+1} = p_j^n
\end{equation}
This system is solved efficiently and unconditionally stably using the Thomas algorithm (tridiagonal matrix algorithm).

At the physical boundaries for velocities (indices \(j=1\) and \(j=N\)), the matrix coefficients are modified so that the flux satisfies:   $F_{1/2} = F_{N+1/2} = 0$. Consequently, the probability is confined within the domain and cannot "leak" through extreme velocity regimes.

Since random disturbances affect the dynamics, the joint mass matrix of the pendulum induces a correlation between the hip and knee accelerations (coefficient \(D_{34}\)). The code accounts for this at the final splitting step using an explicit second-order approximation of the mixed derivative:
\begin{equation}
    p_{i_3, i_4}^{n+1} = p_{i_3, i_4}^n + \Delta t \cdot 2D_{34} \cdot \left[ \frac{p_{i_3+1, i_4+1} - p_{i_3+1, i_4-1} - p_{i_3-1, i_4+1} + p_{i_3-1, i_4-1}}{4\Delta\omega_h \Delta\omega_k} \right]
\end{equation}

After all splitting steps, the mass integral may fluctuate slightly due to round-off errors and the explicit cross-diffusion step. The code enforces a strict normalization at each step:

\begin{equation}
    M = \sum P_{i_1, i_2, i_3, i_4} \cdot (\Delta \theta_h \Delta \theta_k \Delta \omega_h \Delta \omega_k), \quad P \leftarrow \frac{P}{M}
\end{equation}
This guarantees that the total probability is always strictly equal to 1.0 (mass conservation law).

\subsection{Results for a two-link pendulum}

The solution $p(\theta_h, \theta_k, \omega_h, \omega_k, t)$ represents a 4D probability density hypervolume at each time step $t$. To analyze the angular coordinates (\(\theta _{h}\) for the hip joint and \(\theta _{k}\) for the knee joint), the 4D density must be projected (marginalized) onto the angular space by integrating it over the angular velocity grid \(\omega_h, \omega_k\):

\begin{equation}
    p_\theta(\theta_h, \theta_k, t) = \int\limits_{\omega_{h,\min}}^{\omega_{h,\max}} \int\limits_{\omega_{k,\min}}^{\omega_{k,\max}} p(\theta_h, \theta_k, \omega_h, \omega_k, t) \, d\omega_h \, d\omega_k
\end{equation}

The evolution of \(p_{\theta }\) for a two-link pendulum with a small white noise is presented in Figure 4.

\begin{figure}[htbp]
    \centering
     \includegraphics[width=0.8\textwidth]{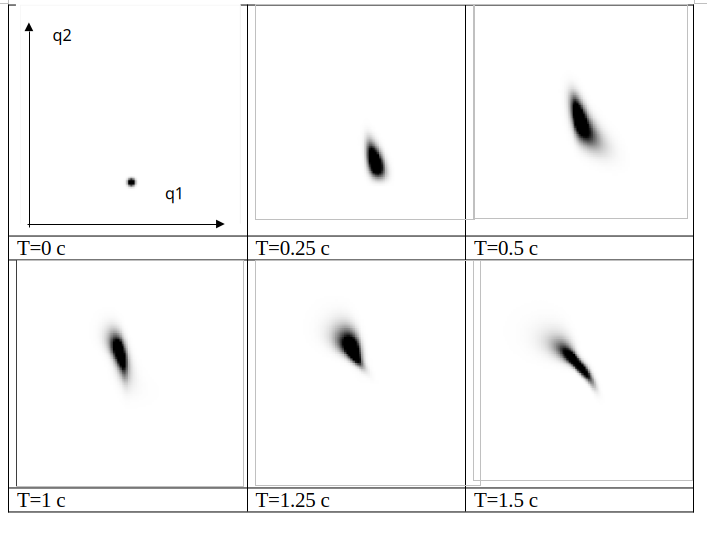} 
    \caption{Probability distribution in the angular space at different time instants. Solution of the boundary value problem.}
    \label{fig:prob_dist_time}
\end{figure}

To obtain the one-dimensional probability density distribution for a specific joint, integration is performed again over the angle of the second joint:
\begin{equation}
    p(\theta_k, t) = \int_{-\pi}^{\pi} p_\theta(\theta_h, \theta_k, t) \, d\theta_h \approx \sum_{i_1=1}^{N_{\theta_h}} P_{i_1, i_2}^{\text{marg}}(t) \cdot \Delta \theta_h
\end{equation}
The evolution of \(p(\theta_k, t)\) is presented in Figure 5. The calculation of \(p(\theta_h, t)\) is carried out analogously.

\begin{figure}[htbp]
    \centering
     \includegraphics[width=1\textwidth]{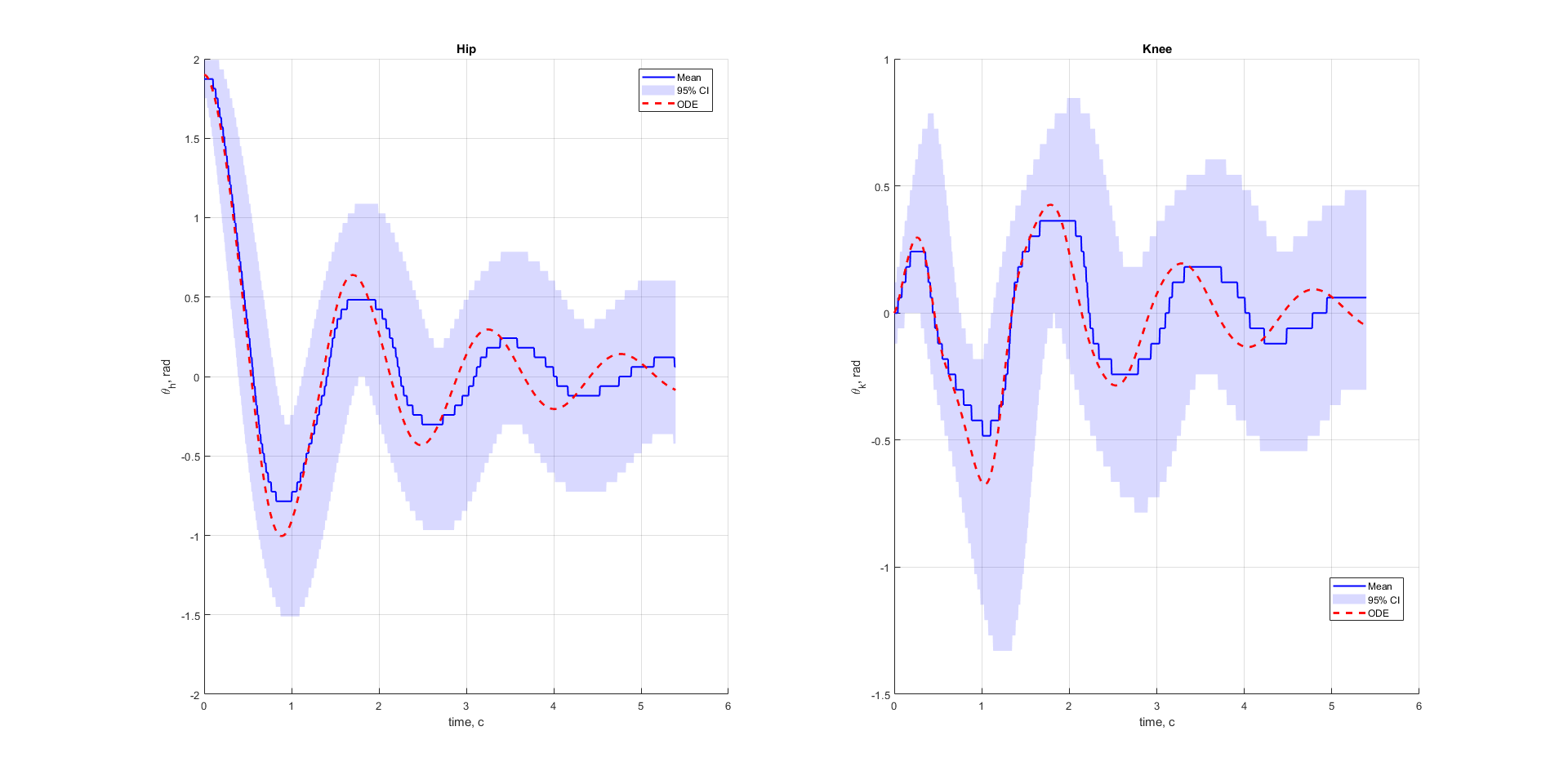} 
    \caption{Comparison of the noiseless system solution and the median value of the probability distribution with a confidence interval for the angular coordinates of a two-link pendulum.}
    \label{fig:double_pendulum_comp}
\end{figure}

\subsection{Results for a gait}

Similarly, one can present the solutions to the Fokker–Planck equation reconstructed from the gait data using the method described above. The solutions are shown in Figures 6 and 7:

\begin{figure}[htbp]
    \centering
     \includegraphics[width=0.8\textwidth]{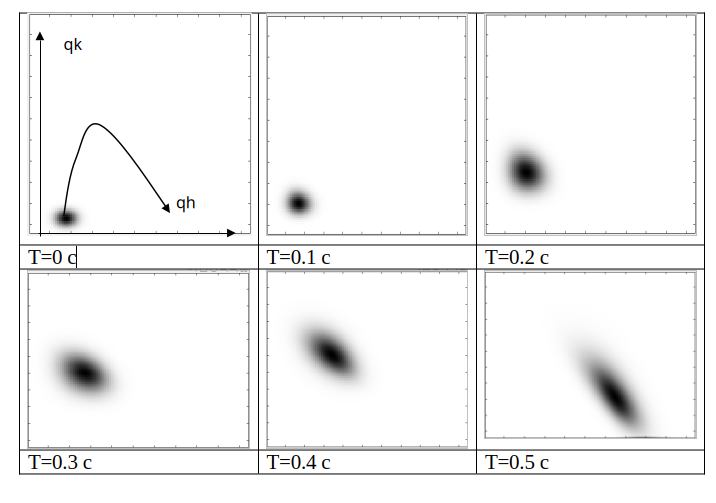}
    \caption{Probability distribution in the gait angular space at different time instants. Solution of the boundary value problem.}
    \label{fig:gait_prob_dist}
\end{figure}

\begin{figure}[htbp]
    \centering
     \includegraphics[width=0.7\textwidth]{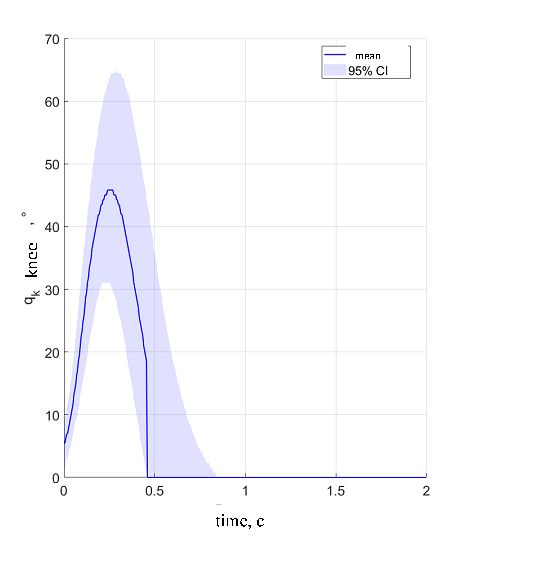}
    \caption{Evolution of the marginal distribution $p(\theta_k,t)$.}
    \label{fig:marginal_evol_k}
\end{figure}

For the time instant corresponding to the maximum of the median value, we estimate the probability distribution. We define the limit state as a probability of no more than 5\% that the maximum knee joint angle falls outside the range of 45–65. The resulting distribution is presented in Figure 8.

\begin{figure}[htbp]
    \centering
     \includegraphics[width=0.8\textwidth]{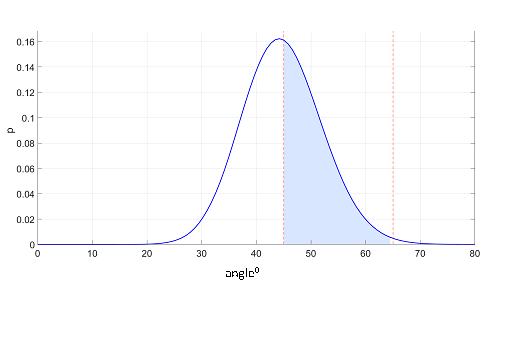}
    \caption{Probability distribution of the maximum knee joint angle realizations at a given spool position.}
    \label{fig:max_knee_angle_dist}
\end{figure}

As can be seen from Figure 8, the limit state is reached:
\begin{equation}
    P\left( (\theta_k < 45) \cup (\theta_k > 65) \right) > 0.05
\end{equation}
This type of CAE analysis serves as an important and highly demanded tool for gait biomechanics investigation and risk-oriented control system synthesis.

\section{Conclusion}

In this paper, a continuous stochastic model of knee joint dynamics is proposed and investigated based on the Itô calculus framework and the Fokker–Planck equation. In contrast to classical approaches, combining the physical Lagrange equations with local nonparametric estimation from big data has enabled the creation of a high-fidelity predictive model that retains physical interpretability.

The introduction of the local Péclet number allowed for a quantitative evaluation of the applicability limits of deterministic prediction in various gait phases. The modal spectral analysis of the Kolmogorov operator, combined with Weyl's asymptotic law, provided a rigorous justification for the computational complexity requirements of prospective neuromorphic and tabular controllers embedded in exoskeleton microcontroller units. The obtained results pave a direct way toward the practical synthesis of stochastic optimal controllers based on the numerical solution of Hamilton–Jacobi–Bellman equations.

\renewcommand{\refname}{Bibliography}

\end{document}